# Brick-wall lattice paths and applications[*]

L. Dăuş[1], V. Beiu[2], S. Cowell[2] and P. Poullin [3]

[1]Department of Mathematics, Technical University of Civil Engineering,
Bdul. Lacul Tei 124, 020396 Bucharest 2, Romania

[2] Department of Mathematics and Computer Science, "Aurel Vlaicu" University,
Str. Elena Dragoi 2-4, 310330 Arad, Romania

[3] Mathematics Division, American University in Dubai,
P.O. Box 28282, Dubai, United Arab Emirates

**Abstract**

We introduce a new type of lattice path, called brick-wall lattice path, and we derive a formula which counts the number of paths on these lattices imposing certain restrictions on the Cartesian plane. Connections to the Fibonacci sequence, as well as to other sequences of numbers, are given. Finally, we use these brick-wall lattice paths to determine the first two non-zero coefficients of the reliability polynomials associated with particular two-terminal networks known as hammocks.

## 1 Introduction

Consider $\mathbf{Z}$ the set of integers and let $S$ be a subset of $\mathbf{Z} \times \mathbf{Z}$. A lattice path with steps in $S$ is a sequence $v_0, v_1, ..., v_k \in \mathbf{Z} \times \mathbf{Z}$ such that each consecutive difference $v_i - v_{i-1}$ lies in $S$.

The first lattice paths enumeration problem can be traced back to the last decades of the 19th century, when the two-candidate ballot problem was proposed by Whitworth [17] in 1878, rediscovered by Bertrand [3] in 1887 and solved by André [1] in the same year. Since then, lattice paths were extensively studied in many fields of mathematics, computer science, and physics [7], [10], and [15].

Quite a few different lattices and various ways of stepping (walking on them) have been proposed, and many of the resultant lattice paths have been analyzed. Still, when looking for lattice paths on a very common type of lattice (which we can see every day around us), we were surprised that such lattice paths have not yet been investigated. The particular lattice we have in mind is the most familiar brick-wall pattern. Although brickwork has generated a wide variety of brick-wall patterns (the interested reader should see https://en.wikipedia.org/wiki/Brickwork), what we are going to analyze in this paper is the pattern/lattice which can be seen in Fig. 1 (technically known as stretcher or running bond). Obviously, the many possible brickwork patterns point to a plethora of lattices on which a large number of lattice paths could be defined and characterized.

[*]This research was funded by the European Union through the European Regional Development Fund under the Competitiveness Operational Program (BioCell-NanoART = Novel Bio-inspired Cellular Nano-architectures, POC-A1-A1.1.4-E nr. 30/2016).

Thus motivated, we have decided to structure this paper as follows. In Section 2 we introduce a new type of lattice path, which, due to its appearance, we call a brick-wall lattice path. Afterwards, we determine matrix formulas counting the number of brick-wall lattice paths bounded by a rectangular region of the Cartesian plane. In Section 3 we present a one-to-one correspondence linking the enumeration problem of brick-wall lattice paths when the bounding rectangle has width 3, on one hand, to the Fibonacci sequence, on the other hand. The proofs are based on a novel $3 \times 3-$matrix of Fibonacci numbers. In Section 4 we use generating functions to determine the number of brick-wall lattice paths when the width of the bounding rectangle is 4. Section 5 indicates several entries in the OEIS [18] that are counting sequences for certain combinatorial classes of brick-wall lattice paths. Finally, Section 6 reveals links between the results detailed above and the reliability polynomials associated with hammock networks (which have a brick-wall representation). Although this particular type of network was proposed and studied more than sixty years ago by Moore and Shannon [11], the associated reliability polynomials were expressed exactly only for a few particular cases [6]. Using the counting process determined in Section 2, we will prove exact formulas for the first and second non-zero coefficients of the reliability polynomials for any hammock networks.

## 2 Brick-wall lattice paths

**Definition 2.1.** *We say that the lattice point $(a, b) \in \mathbf{Z} \times \mathbf{Z}$ is odd (even) if $|a| + |b|$ is an odd (resp. even) number.*

**Definition 2.2.** *A brick-wall lattice path of the first (second) type is a lattice path that consists of the following steps: horizontal $H = (1, 0)$, up $U_e = (0, 1)$ for even points and down $D_o = (0, -1)$ for odd points (resp. horizontal $H = (1, 0)$, up $U_o = (0, 1)$ for odd points and down $D_e = (0, -1)$ for even points).*

Drawing all the possible brick-wall lattice paths of the first type (or of the second type), will "fill" the Cartesian plane with non-overlapping horizontal bricks having length 2 (two steps on the horizontal direction) and width 1 (one step on the vertical direction). The first type of brick-wall lattice path corresponds to the case when the origin is in the corner of a brick, while the second type matches the case when the origin is the midpoint on the length of a brick. Figure 1 illustrates two types of brick-wall lattice paths from $(0, 3)$ to $(12, 7)$.

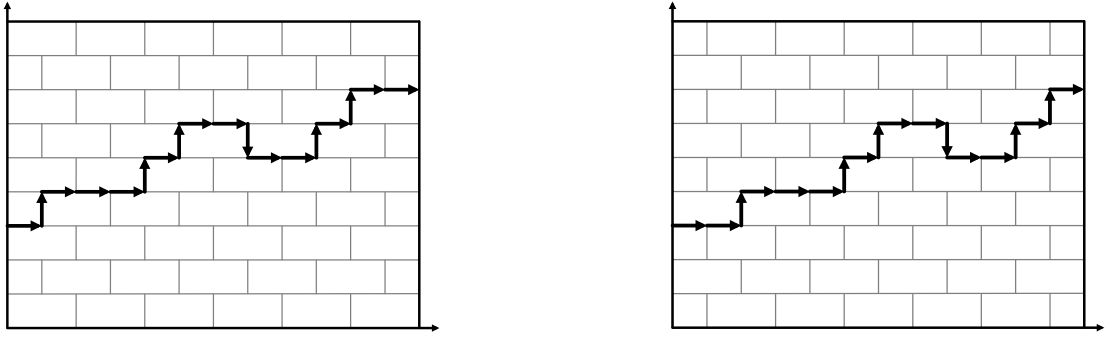

| a) Sample of brick-wall lattice path of the first type | b) Sample of brick-wall lattice path of the second type |

Figure 1: Two different brick-wall lattice paths

Let $l, w$ (i.e., length and width) be two integers, $l, w \geq 2$, and consider the rectangle $\mathcal{R} = [0, l-1] \times [0, w-1]$. We denote by $\mathcal{B}^{(1)}_{l,w}$ (and $\mathcal{B}^{(2)}_{l,w}$) the number of all brick-wall lattice paths of the first (resp. second) type connecting a point $(0, \alpha), \alpha \in \{0, 1, ..., w-1\}$ (the left side of $\mathcal{R}$) with a point $(l-1, \beta), \beta \in \{0, 1, ..., w-1\}$ (the right side of $\mathcal{R}$) with walks inside $\mathcal{R}$. Therefore, all the brick-wall lattice paths we are considering here are restricted to the rectangle $\mathcal{R}$. Moreover, in the paths we are counting, no vertical steps are allowed along the left side of $\mathcal{R}$.

For any $0 \leq i \leq l-1$ and any $0 \leq j \leq w-1$, we denote by $b^{(1)}_{ij}$ (and $b^{(2)}_{ij}$) the number of brick-wall lattice paths of the first (resp. second) type from the left side of $\mathcal{R}$ to $(i, j)$. We set

$$b^{(t)}_i = \begin{pmatrix} b^{(t)}_{i,w-1} \\ b^{(t)}_{i,w-2} \\ \cdots \\ b^{(t)}_{i,2} \\ b^{(t)}_{i,1} \\ b^{(t)}_{i,0} \end{pmatrix},$$

for any $0 \leq i \leq l-1$ and $t \in \{1, 2\}$. Considering the trivial (empty) paths,

$$b^{(1)}_0 = b^{(2)}_0 = \begin{pmatrix} 1 \\ 1 \\ \cdots \\ 1 \\ 1 \\ 1 \end{pmatrix}.$$

We denote $b^{(1)}_0 = b^{(2)}_0$ by $b_0$.

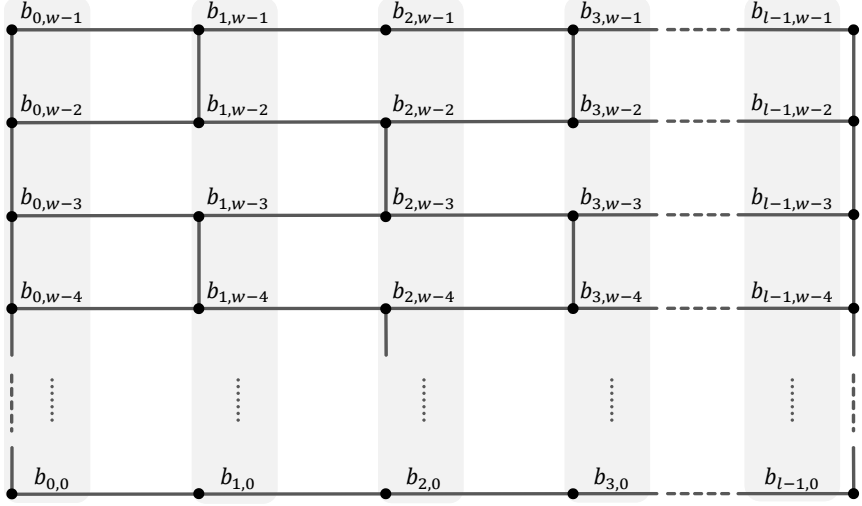

Figure 2: Brick-wall lattice

## 2.1 If $w$ is an odd number

In this case, any brick-wall lattice path of the first type corresponds to one and only one brick-wall lattice path of the second type, the two paths being symmetric one each other with respect the horizontal line $y = \frac{w-1}{2}$. Therefore $\mathcal{B}_{l,w}^{(1)} = \mathcal{B}_{l,w}^{(2)}$. Let us denote this common value by $\mathcal{B}_{l,w}$. To compute $\mathcal{B}_{l,w}$ we will use two particular endomorphisms of the $w-$dimensional real vector space $\mathbf{R}^w$.

The map from numbers of the column $b_{2k}^{(1)}$ to numbers of the next column $b_{2k+1}^{(1)}$ is the endomorphism

$$T_U : \mathbf{R}^w \longrightarrow \mathbf{R}^w,$$

$$T_U(x_1, x_2, x_3, x_4, ..., x_{w-2}, x_{w-1}, x_w) = \\ (x_1 + x_2, x_1 + x_2, x_3 + x_4, x_3 + x_4, ..., x_{w-2} + x_{w-1}, x_w),$$

represented by the $w \times w$-matrix:

$$M_U = \begin{pmatrix} 1 & 1 & 0 & 0 & ... & 0 & 0 & 0 \\ 1 & 1 & 0 & 0 & ... & 0 & 0 & 0 \\ 0 & 0 & 1 & 1 & ... & 0 & 0 & 0 \\ 0 & 0 & 1 & 1 & ... & 0 & 0 & 0 \\ ... & ... & ... & ... & ... & ... & ... & ... \\ 0 & 0 & 0 & 0 & ... & 1 & 1 & 0 \\ 0 & 0 & 0 & 0 & ... & 1 & 1 & 0 \\ 0 & 0 & 0 & 0 & ... & 0 & 0 & 1 \end{pmatrix}. \qquad (1)$$

Therefore
$$M_U b^{(1)}_{2k} = b^{(1)}_{2k+1}, \qquad (2)$$
for any $0 \leq k \leq \left\lfloor \dfrac{l-2}{2} \right\rfloor$.

Similarly, to evaluate the values of the column $b^{(1)}_{2k+2}$ from the values of the previous column $b^{(1)}_{2k+1}$ we use the endomorphism

$$T_L : \mathbf{R}^w \longrightarrow \mathbf{R}^w,$$

$$T_L (x_1, x_2, x_3, ..., x_{w-1}, x_w) = (x_1, x_2 + x_3, x_2 + x_3, x_4 + x_5 x_4 + x_5, ..., x_{w-1} + x_w),$$

given by the $w \times w$-matrix:

$$M_L = \begin{pmatrix} 1 & 0 & 0 & 0 & 0 & ... & 0 & 0 \\ 0 & 1 & 1 & 0 & 0 & ... & 0 & 0 \\ 0 & 1 & 1 & 0 & 0 & ... & 0 & 0 \\ 0 & 0 & 0 & 1 & 1 & ... & 0 & 0 \\ 0 & 0 & 0 & 1 & 1 & ... & 0 & 0 \\ ... & ... & ... & ... & ... & ... & ... & ... \\ 0 & 0 & 0 & 0 & 0 & ... & 1 & 1 \\ 0 & 0 & 0 & 0 & 0 & ... & 1 & 1 \end{pmatrix}. \qquad (3)$$

Hence
$$M_L b^{(1)}_{2k+1} = b^{(1)}_{2k+2}, \qquad (4)$$
for any $0 \leq k \leq \left\lfloor \dfrac{l-3}{2} \right\rfloor$.

Using (2), (4), together with an inductive argument, we immediately obtain the following result:

**Theorem 2.3.** *When $w$ is an odd number, then for any $1 \leq r \leq l-1$ the components of the column $b^{(1)}_r$ are given by*

$$b^{(1)}_r = \begin{cases} (M_L M_U)^\rho b_0 & , \text{ if } r \text{ is even} \\ M_U (M_L M_U)^\rho b_0 & , \text{ if } r \text{ is odd} \end{cases},$$

*where $\rho = \left\lfloor \dfrac{r}{2} \right\rfloor$. Hence the number of all brick-wall lattice paths of the first (as well as second) type connecting the left and right sides of the rectangle $\mathcal{R}$, with walks inside $\mathcal{R}$, is*

$$\mathcal{B}_{l,w} = u b^{(1)}_{l-1}, \qquad (5)$$

*where $u = (1, 1, ..., 1)$ is a row matrix having $w$ components.*

## 2.2 If $w$ is even number

When $w$ is an even number, but $l$ is odd, then the brick-wall lattice (on which we will count the paths) of the first type exhibits $\left\lfloor \dfrac{(w-1)(l-2)}{2} \right\rfloor$ possible vertical steps, while the brick-wall lattice of the second type has one more: $\left\lceil \dfrac{(w-1)(l-2)}{2} \right\rceil$. Therefore $\mathcal{B}_{l,w}^{(1)} \neq \mathcal{B}_{l,w}^{(2)}$. Hence, when $w$ is even we need to study the two types of brick-wall lattices separately. Fig. 3 shows the two types of brick-wall lattice paths for the case $w = 4$ and $l = 5$: the first one with 4 possible vertical steps, while the second one having 5.

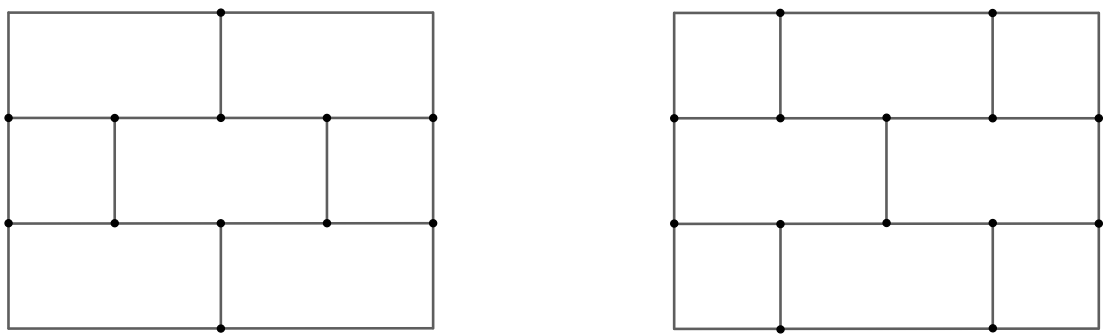

Figure 3: Brick-wall lattice paths of the first and second type for w=4 and l=5

### 2.2.1 Brick-wall lattice paths of the first type

We may go from the left side to the right side of the rectangle $\mathcal{R}$ by using, alternatively, the linear transformation given by the $w \times w$-matrix:

$$M_- = \begin{pmatrix} 1 & 0 & 0 & \ldots & 0 & 0 & 0 \\ 0 & 1 & 1 & \ldots & 0 & 0 & 0 \\ 0 & 1 & 1 & \ldots & 0 & 0 & 0 \\ \ldots & \ldots & \ldots & \ldots & \ldots & \ldots & \ldots \\ 0 & 0 & 0 & \ldots & 1 & 1 & 0 \\ 0 & 0 & 0 & \ldots & 1 & 1 & 0 \\ 0 & 0 & 0 & \ldots & 0 & 0 & 1 \end{pmatrix}, \qquad (6)$$

followed by another linear transformation corresponding to the $w \times w$-matrix:

$$M_+ = \begin{pmatrix} 1 & 1 & 0 & 0 & \ldots & 0 & 0 & 0 & 0 \\ 1 & 1 & 0 & 0 & \ldots & 0 & 0 & 0 & 0 \\ 0 & 0 & 1 & 1 & \ldots & 0 & 0 & 0 & 0 \\ 0 & 0 & 1 & 1 & \ldots & 0 & 0 & 0 & 0 \\ \ldots & \ldots & \ldots & \ldots & \ldots & \ldots & \ldots & \ldots & \ldots \\ 0 & 0 & 0 & 0 & \ldots & 1 & 1 & 0 & 0 \\ 0 & 0 & 0 & 0 & \ldots & 1 & 1 & 0 & 0 \\ 0 & 0 & 0 & 0 & \ldots & 0 & 0 & 1 & 1 \\ 0 & 0 & 0 & 0 & \ldots & 0 & 0 & 1 & 1 \end{pmatrix}. \qquad (7)$$

Obviously, for any $1 \leq r \leq l-1$ the components of the column $b_r^{(1)}$ are given by

$$b_r^{(1)} = \begin{cases} (M_+M_-)^\rho b_0 & \text{, if } r \text{ is even} \\ M_-(M_+M_-)^\rho b_0 & \text{, if } r \text{ is odd} \end{cases}, \qquad (8)$$

where $\rho = \left\lfloor \dfrac{r}{2} \right\rfloor$.

### 2.2.2 Brick-wall lattice paths of the second type

In order to count this type of brick-wall lattice path it is necessary to use, alternatively, the matrices $M_+$ and $M_-$, introduced in (7) and (6), respectively. Therefore, to obtain the components of the column $b_r^{(2)}$ we have to modify (8) by interchanging the two matrices $M_+$ and $M_-$.

We can summarize what we have shown so far in:

**Theorem 2.4.** *Let $w$ be an even number, $u = (1,1,...,1) \in \mathbf{R}^w$ and $\rho = \left\lfloor \dfrac{l-1}{2} \right\rfloor$. Then:*

*i) the number of brick-wall lattice paths of the first type is*

$$\mathcal{B}_{l,w}^{(1)} = ub_{l-1}^{(1)}, \qquad (9)$$

*where*

$$b_{l-1}^{(1)} = \begin{cases} (M_+M_-)^\rho b_0 & \text{, if } l \text{ is odd} \\ M_-(M_+M_-)^\rho b_0 & \text{, if } l \text{ is even} \end{cases};$$

*ii) the number of brick-wall lattice paths of the second type is*

$$\mathcal{B}_{l,w}^{(2)} = ub_{l-1}^{(2)}, \qquad (10)$$

*where*

$$b_{l-1}^{(2)} = \begin{cases} (M_-M_+)^\rho b_0 & \text{, if } l \text{ is odd} \\ M_+(M_-M_+)^\rho b_0 & \text{, if } l \text{ is even} \end{cases}.$$

## 3 Connections with Fibonacci numbers

In this section we investigate the sequence $(\mathcal{B}_{l,3})_{l \geq 0}$ which counts both types of brick-wall lattice paths related to rectangles $\mathcal{R}$ of width $w = 3$. For the particular $3 \times 3$-matrices $M_L$ and $M_U$ given by (3) and (1), respectively, we obtain

$$M_L M_U = \begin{pmatrix} 1 & 0 & 0 \\ 0 & 1 & 1 \\ 0 & 1 & 1 \end{pmatrix} \begin{pmatrix} 1 & 1 & 0 \\ 1 & 1 & 0 \\ 0 & 0 & 1 \end{pmatrix} = \begin{pmatrix} 1 & 1 & 0 \\ 1 & 1 & 1 \\ 1 & 1 & 1 \end{pmatrix}. \qquad (11)$$

Using the Fibonacci sequence $F_0 = 0, F_1 = 1, F_n = F_{n-1} + F_{n-2}$, for any $n \geq 2$, we may rewrite the product $M_L M_U$ in (11) as

$$\begin{pmatrix} 1 & 1 & 0 \\ 1 & 1 & 1 \\ 1 & 1 & 1 \end{pmatrix} = \begin{pmatrix} F_1 & F_1 & F_0 \\ F_2 & F_2 & F_1 \\ F_2 & F_2 & F_1 \end{pmatrix}. \tag{12}$$

If we denote the previous matrix by $\mathcal{F}$, a direct computation gives:

$$\mathcal{F}^2 = \begin{pmatrix} 2 & 2 & 1 \\ 3 & 3 & 2 \\ 3 & 3 & 2 \end{pmatrix} = \begin{pmatrix} F_3 & F_3 & F_2 \\ F_4 & F_4 & F_3 \\ F_4 & F_4 & F_3 \end{pmatrix}$$

and

$$\mathcal{F}^3 = \begin{pmatrix} 5 & 5 & 3 \\ 8 & 8 & 5 \\ 8 & 8 & 5 \end{pmatrix} = \begin{pmatrix} F_5 & F_5 & F_4 \\ F_6 & F_6 & F_5 \\ F_6 & F_6 & F_5 \end{pmatrix}.$$

Following on these examples, we prove the next results.

**Proposition 3.1.** *Using the previous notations, for any integer $n \geq 1$ the following relations hold:*

$$\mathcal{F}^n = \begin{pmatrix} F_{2n-1} & F_{2n-1} & F_{2n-2} \\ F_{2n} & F_{2n} & F_{2n-1} \\ F_{2n} & F_{2n} & F_{2n-1} \end{pmatrix} \tag{13}$$

*and*

$$M_U \mathcal{F}^n = \begin{pmatrix} F_{2n+1} & F_{2n+1} & F_{2n} \\ F_{2n+1} & F_{2n+1} & F_{2n} \\ F_{2n} & F_{2n} & F_{2n-1} \end{pmatrix}. \tag{14}$$

*Proof.* By induction on $n$. □

**Corollary 3.2.** *For any $l \geq 1$, the following identity holds:*

$$\mathcal{B}_{l,3} = F_{l+3}. \tag{15}$$

*Proof.* For $l = 1$, $\mathcal{R}$ is a degenerate rectangle (a vertical line segment), hence $\mathcal{B}_{1,3} = 3 = F_{1+3}$. When $l \geq 2$, we first consider $l$ as an odd number. Then $l = 2k+1$ and $\left\lfloor \dfrac{l-1}{2} \right\rfloor = k$. By *Theorem 1* and *Proposition 1*, the number of the brick-wall lattice paths is

$$u \mathcal{F}^k b_0 = \begin{pmatrix} 1 & 1 & 1 \end{pmatrix} \begin{pmatrix} F_{2k-1} & F_{2k-1} & F_{2k-2} \\ F_{2k} & F_{2k} & F_{2k-1} \\ F_{2k} & F_{2k} & F_{2k-1} \end{pmatrix} \begin{pmatrix} 1 \\ 1 \\ 1 \end{pmatrix}$$

$$= 4F_{2k} + 4F_{2k-1} + F_{2k-2} = F_{2k+4} = F_{l+3}.$$

If $l$ is an even number, $l = 2k$, then $\left\lfloor \frac{l-1}{2} \right\rfloor = k - 1$. In this case, applying again *Theorem 1* and *Proposition 1*, the number of the brick wall paths is

$$uM_U \mathcal{F}^{k-1} b_0 = \begin{pmatrix} 1 & 1 & 1 \end{pmatrix} \begin{pmatrix} F_{2k-1} & F_{2k-1} & F_{2k-2} \\ F_{2k-1} & F_{2k-1} & F_{2k-2} \\ F_{2k-2} & F_{2k-2} & F_{2k-3} \end{pmatrix} \begin{pmatrix} 1 \\ 1 \\ 1 \end{pmatrix}$$

$$= 4F_{2k-1} + 4F_{2k-2} + F_{2k-3} = F_{2k+3} = F_{l+3}.$$

□

*Remarks.* In 1951 Brenner [4] used for the first time a matrix related to the Fibonacci sequence. That matrix is now known as Fibonacci $Q$-matrix and it has the following entries:

$$Q = \begin{pmatrix} 1 & 1 \\ 1 & 0 \end{pmatrix}. \tag{16}$$

Obviously, the matrix $Q$ is unimodular and its eigenvalues are $\varphi$ and $-\varphi^{-1}$. The matrix $\mathcal{F}$ introduced in (12) is singular and it has $0, \varphi^2$ and $\varphi^{-2}$ as eigenvalues.

As is well known, the Fibonacci sequence $(F_n)_{n=0}^{\infty}$ has generating function $F(z) = \frac{z}{1-z-z^2}$. That is to say, for all integers $n \geq 0$,

$$F_n = [z^n] \frac{z}{1 - z - z^2}$$

where $[z^n]F(z)$ denotes the coefficient of $z^n$ in the formal power series expansion of $F(z)$ centred at $z = 0$. Therefore, the number of brick-wall lattice paths in the case $w = 3$ is

$$F_{l+3} = [z^{l+3}]F(z) = [z^l] \frac{F(z)}{z^3} = [z^l] \frac{1}{z^2 - z^3 - z^4}.$$

## 4 Brick-wall paths of width $w = 4$

In this section, and also in the next one, by counting the brick-wall lattice paths for different values of $w$, we will start to recover several sequences that appear in *Sloane's On-Line Encyclopedia of Integer Sequences* [18].

Let us consider first the brick-wall lattice paths of the second type restricted to the rectangle $\mathcal{R}$. For all of them, there are $\left\lceil \frac{(l-2)(w-1)}{2} \right\rceil = \left\lceil \frac{3(l-2)}{2} \right\rceil$ possible vertical steps. In this case, the two matrices $M_-$ and $M_+$ are as described in (6) and (7), respectively. Moreover, $M_- M_+$ is diagonalisable, with eigenvalues $(0, 0, 1, 3)$. It follows that the number of brick-wall lattice paths $\mathcal{B}_{l,4}^{(2)}$ is given by

$$\mathcal{B}_{2k+1,4}^{(2)} = \alpha_1 1^k + \beta_1 3^k = \alpha_1 + \beta_1 3^k,$$

when $l = 2k + 1$ is odd, and $\mathcal{B}_{l,4}^{(2)}$ is given by

$$\mathcal{B}_{2k,4}^{(2)} = \alpha_2 1^k + \beta_2 3^k = \alpha_2 + \beta_2 3^k,$$

when $l = 2k$ is even. Here the $\alpha_i$ and $\beta_i$ are constants to be determined. In fact, since

$$\mathcal{B}^{(2)}_{1,4} = \begin{pmatrix} 1 & 1 & 1 & 1 \end{pmatrix} b_0 = \begin{pmatrix} 1 & 1 & 1 & 1 \end{pmatrix} \begin{pmatrix} 1 \\ 1 \\ 1 \\ 1 \end{pmatrix} = 4$$

and

$$\mathcal{B}^{(2)}_{3,4} = \begin{pmatrix} 1 & 1 & 1 & 1 \end{pmatrix} b_2 = \begin{pmatrix} 1 & 1 & 1 & 1 \end{pmatrix} \begin{pmatrix} 2 \\ 4 \\ 4 \\ 2 \end{pmatrix} = 12$$

we obtain $\alpha_1 = 0$ and $\beta_1 = 4$, hence

$$\mathcal{B}^{(2)}_{2k+1,4} = 4 \cdot 3^k, \qquad k \in \{0, 1, 2, \ldots\}.$$

Similarly, since $\mathcal{B}^{(2)}_{2,4} = 8$ and $\mathcal{B}^{(2)}_{4,4} = 24$ we obtain $\alpha_2 = 0$ and $\beta_2 = \frac{8}{3}$, hence

$$\mathcal{B}^{(2)}_{2k,4} = \frac{8}{3} \cdot 3^k = 8 \cdot 3^{k-1}, \qquad k \in \{1, 2, 3, \ldots\}.$$

The sequence of numbers of brick-wall lattice paths for odd values of $l$ therefore has generating function

$$G^{(2)}_{4,\text{odd}}(z) = \sum_{k=0}^{\infty} \mathcal{B}^{(2)}_{2k+1,4} z^k = \sum_{k=0}^{\infty} 4 \cdot 3^k z^k = \frac{4}{1 - 3z},$$

and the corresponding generating function for even values of $l$ is

$$G^{(2)}_{4,\text{even}}(z) = \sum_{k=1}^{\infty} \mathcal{B}^{(2)}_{2k,4} z^k = \sum_{k=1}^{\infty} 8 \cdot 3^{k-1} z^k = \frac{8z}{1 - 3z}.$$

Therefore the generating function for the sequence of numbers of brick-wall lattice paths for arbitrary $l$ is

$$G^{(2)}_4(z) = \sum_{l=1}^{\infty} \mathcal{B}^{(2)}_{l,4} z^l = \sum_{k=0}^{\infty} \mathcal{B}^{(2)}_{2k+1,4} z^{2k+1} + \sum_{k=1}^{\infty} \mathcal{B}^{(2)}_{2k,4} z^{2k}$$

$$= z \sum_{k=0}^{\infty} \mathcal{B}^{(2)}_{2k+1,4} (z^2)^k + \sum_{k=1}^{\infty} \mathcal{B}^{(2)}_{2k,4} (z^2)^k = z \cdot G^{(2)}_{4,\text{odd}}(z^2) + G^{(2)}_{4,\text{even}}(z^2)$$

$$= z \frac{4}{1 - 3z^2} + \frac{8z^2}{1 - 3z^2} = \frac{4z + 8z^2}{1 - 3z^2}.$$

For brick-wall lattice paths of the first type restricted to $\mathcal{R}$, there are $\left\lfloor \frac{(l-2)(w-1)}{2} \right\rfloor = \left\lfloor \frac{3(l-2)}{2} \right\rfloor$ possible vertical steps. In this alternate case, $M_+$ and $M_-$ must be interchanged, and then a completely similar analysis shows that the corresponding generating function is

$$G^{(1)}_4(z) = \frac{4z + 6z^2}{1 - 3z^2}.$$

Thus the numbers of brick-wall lattice paths are the coefficients appearing in the right hand side of

$$G_4^{(2)}(z) = \frac{4z + 8z^2}{1 - 3z^2} = 4z + 8z^2 + 12z^3 + 24z^4 + 36z^5 + 72z^6 + 108z^7 + 216z^8 + \cdots \quad (17)$$

in the first case, and those coefficients appearing in the right hand side of

$$G_4^{(1)}(z) = \frac{4z + 6z^2}{1 - 3z^2} = 4z + 6z^2 + 12z^3 + 18z^4 + 36z^5 + 54z^6 + 108z^7 + 162z^8 + \cdots \quad (18)$$

in the alternate case.

Note that if both $w$ and $l$ are even, then by symmetry, the path counts are the same for both types.

**Remark 4.1.** *The sequence of numbers of brick wall paths given by $G_4^+$ in (15) is the sequence [A153339] [18] which is the number of zig-zag paths from top to bottom of a rectangle of width 5 with n rows whose color is that of the top right corner.*

*The sequence of numbers of brick wall paths given by $G_4$ in (16) is the sequence [A068911] [18] which can be described as the number of n step walks (each step ±1 starting from 0) which are never more than 2 or less than -2 (see also [13]).*

## 5 Other links to OEIS

In Section 3 we have seen that the number of brick-wall lattice paths for the case when $w = 3$ is $F_{l+4}$ [14], [16], which is [A000045] in OEIS [18]. It was a normal step to check if the numbers of brick-wall lattice paths for the case when $w = 4$ are also in OEIS. The answer was positive, both for $\mathcal{B}_{l,4}^{(1)}$ and for $\mathcal{B}_{l,4}^{(2)}$. That is why we have decided to look in OEIS for the numbers of brick-wall lattice paths for rectangles $\mathcal{R}$ of other widths $w$. We were able to find sequences in OEIS matching the numbers of brick-wall lattice paths for $w = 3, 4, 5, \ldots, 11$. All of these are presented in a compact form in Table 1 where we have included not only the number of the series in OEIS, but, for completeness, a few references also.

Last but not least, the first five sequences for odd $w$, i.e., $w = 3, 5, 7, 9, 11$, form an array which is matching the even rows $k = 0, 2, 4, 6, 8$ of [A116183] (for $k = 8$ it is mentioned that the sequence is [A129638]). As [A116183] is not defined beyond $k = 9$ (as sequence [A129639]) we can only speculate the possibility of an identity. Still, for even $w$ and odd $k$ the sequences are different at least for the fact that for $w$ even we have two different sequences.

It is obvious that, for any fixed $w$, the number of brick-wall lattice paths generate another sequence by varying $l$, so the matrix formulas (5), (9), and (10) determined in Section 2 are generating an infinite number of sequences.

Table 1:

| $\mathcal{B}_{l,w}$ | First eight terms | Sequence in OEIS [18] | Other references |
|---|---|---|---|
| $\mathcal{B}_{l,3}$ | 3, 5, 8, 13, 21, 34, 55, 89 | [A000045] | [14], [16] |
| $\mathcal{B}_{l,4}^{(1)}$ | 4, 6, 12, 18, 36, 54, 108, 162 | [A068911] | [13] |
| $\mathcal{B}_{l,4}^{(2)}$ | 4, 8, 12, 24, 36, 72, 108, 216 | [A153339] | [12] |
| $\mathcal{B}_{l,5}$ | 5, 9, 16, 29, 52, 94, 169, 305 | [A090990] | [8], [9] |
| $\mathcal{B}_{l,6}^{(1)}$ | 6, 10, 20, 34, 68, 116, 232, 396 | [A030436] | [2] |
| $\mathcal{B}_{l,6}^{(2)}$ | 6, 12, 20, 40, 68, 136, 232, 464 | [A030435] | [12] |
| $\mathcal{B}_{l,7}$ | 7, 13, 24, 45, 84, 158, 296, 557 | [A090992] | [8], [9] |
| $\mathcal{B}_{l,8}^{(1)}$ | 8, 14, 28, 50, 100, 180, 360, 650 | [A153364] | [12] |
| $\mathcal{B}_{l,8}^{(2)}$ | 8, 16, 28, 56, 100, 200, 360, 720 | [A153363] | [12] |
| $\mathcal{B}_{l,9}$ | 9, 17, 32, 61, 116, 222, 424, 813 | [A090994] | [8], [9] |
| $\mathcal{B}_{l,10}^{(1)}$ | 10, 18, 36, 66, 132, 244, 488, 906 | [A153370] | [12] |
| $\mathcal{B}_{l,10}^{(2)}$ | 10, 20, 36, 72, 132, 264, 488, 976 | [A153369] | [12] |
| $\mathcal{B}_{l,11}$ | 11, 21, 40, 77, 148, 286, 552, 1069 | [A129638] | [8], [9] |

# 6 Reliability polynomial of brick-wall network

In reliability theory hammock networks were introduced by Moore and Shannon [11]. Their idea was to replace an unreliable two-contact electromechanical relay by a (redundant) network of such relays. Such a network is known as a two-terminal network because it has an input/source $S$, and an output/terminus $T$. They also assumed that all the relays in such networks are switched simultaneously (i.e., are activated, or not activated, simultaneously). The relays forming the network were represented as sequences of independent identically distributed (i.i.d.) random variables, being closed with probability $p$ and open with probability $q = 1-p$. It follows that the state of a two-terminal network is a Boolean function of random variables, and the probability of the network being closed is a function $h(p)$ known as the reliability polynomial associated with the network.

In general, for an arbitrary network, there are several forms for representing the reliability polynomial $h(p)$, but the most used one is known as *N-form* [5], [11]:

$$h(p) = \sum_{k=0}^{m} N_k p^k (1-p)^{m-k},$$

where $m$ is the total number of contacts of the network and $N_k$ is the number of ways we can select a subset of $k$ contacts in the network such that if these $k$ contacts are closed, and the remaining contacts open, then the network is closed. Obviously, if we expand the reliability polynomial in powers of $p$ we obtain the *P-form*:

$$h(p) = \sum_{k=0}^{m} P_k p^k$$

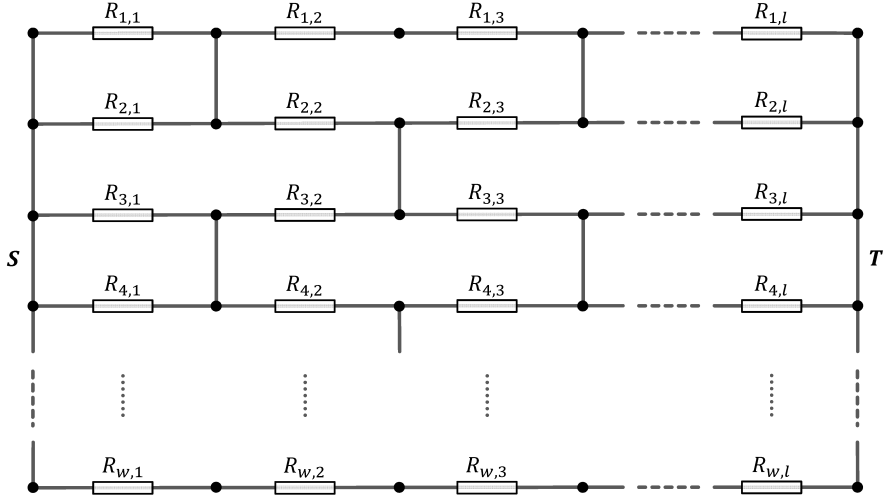

Figure 4: Hammock network

In the particular case of a hammock network of length $l$ and width $w$, the *N-form* of the reliability polynomial may be written as

$$h(p) = \sum_{k=l}^{wl} N_k p^k (1-p)^{wl-k}, \qquad (19)$$

while the *P-form* is

$$h(p) = \sum_{k=l}^{wl} P_k p^k. \qquad (20)$$

Using (19) and (20), the coefficients of $p^l$ and $p^{l+1}$ are then, respectively,

$$\begin{array}{ll} P_l = N_l & (\textit{first coefficient}) \\ P_{l+1} = -(wl-l)N_l + N_{l+1} & (\textit{second coefficient}) \end{array}$$

Since $N_k$ represents the number of $k$-subsets of relays containing a path from $S$ to $T$, then $N_l$ is exactly the number of brick-wall lattice paths inside the rectangle $\mathcal{R}$ having length $l$ and width $w$. Hence

$$N_l = \mathcal{B}_{l,w}^{(1)} \quad \text{or} \quad N_l = \mathcal{B}_{l,w}^{(2)}, \qquad (21)$$

according to the distribution of the vertical connections inside the hammock network.

The $(l+1)$-subsets of relays counted by $N_{l+1}$ are obtained from the selection of a path from $S$ to $T$ ($P_l = N_l$ choices) and of a supplemental relay ($wl - l$ choices). It gives $(wl - l)N_l$ choices minus the overcounting, which occurs when a half-brick at any edge of the wall is fully selected (see Fig. 1). The overcounted $(l+1)$-subsets contain actually two paths from $S$ to $T$, and hence, the number of such $(l+1)$-subsets is half the number of paths from $S$ to $T$ starting (resp. finishing) along a half-brick. Denoting by $M$ the transition matrix computed in *Theorem 1* or *Theorem 2*, such paths are easily counted by adjusting the column vector (resp. row vector) in the expression

$$\begin{pmatrix} 1 & \cdots & 1 \end{pmatrix} M \begin{pmatrix} 1 \\ \vdots \\ 1 \end{pmatrix} = uMu^T$$

so the relays on the left (resp. right) edge of the wall that belong to a half-brick yield an entry 0 instead of 1. There are only a few such relays (0, 1, or 2), located at the top or bottom of the wall. It yields an expression of the form

$$uMu_1^T + u_2 Mu^T,$$

where $u_1$ (resp. $u_2$) is the modified column (resp. row).

Therefore,
$$N_{l+1} = (wl - l)N_l - \frac{1}{2}(uMu_1^T + uM^T u_2^T),$$

and hence
$$P_{l+1} = -\frac{1}{2}(uMu_1^T + uM^T u_2^T).$$

It specializes as follows. Let
$$\begin{aligned} u &= \begin{pmatrix} 1 & 1 & 1 & \cdots & 1 & 1 & 1 \end{pmatrix}, \\ v &= \begin{pmatrix} 1 & 1 & 1 & \cdots & 1 & 1 & 0 \end{pmatrix}, \\ t &= \begin{pmatrix} 0 & 1 & 1 & \cdots & 1 & 1 & 0 \end{pmatrix}. \end{aligned}$$

For odd width, the transition matrix is (wlog) $M = M_L M_U \ldots M_L M_U$ or $M = M_U M_L M_U \ldots M_L M_U$. In both cases,
$$P_{l+1} = -uMv^T. \tag{22}$$

For even width, the transition matrix may be $M = M_- M_+ \ldots M_- M_+ M_-$, in which case
$$P_{l+1} = -uMt^T. \tag{23}$$

It may be $M = M_+ M_- \ldots M_+ M_- M_+$, in which case
$$P_{l+1} = -uMu^T. \tag{24}$$

Finally, it may be (wlog) $M = M_+ M_- \ldots M_+ M_-$, in which case
$$P_{l+1} = -\frac{1}{2}(uMt^T + uM^T u^T). \tag{25}$$

**Remark 6.1.** *Link with Fibonacci: for $w = 3$, a direct computation (similar to the one for the first coefficient) then shows that $P_{l+1} = -2F_{l+1}$ independently of the parity of $l$.*